\documentclass[runningheads,a4paper]{llncs}

\usepackage{amsmath,amssymb,amsbsy,latexsym,mathtools}
\usepackage{graphicx,epic,eepic}
\usepackage{cite}
\usepackage{comment}
\usepackage{wrapfig}
\usepackage{epstopdf}
\usepackage{xpatch}
\usepackage[intoc]{nomencl}

\usepackage{etoolbox,ragged2e,siunitx}

\renewcommand\nomgroup[1]{\def\nomtemp{\csname nomstart#1\endcsname}\nomtemp}

\newcommand{\nomDef}[1]{\parbox[t]{10cm}{\RaggedRight #1}}

\newcommand{\nomtypeD}[3][]{\nomenclature[D#1]{#2}{\nomDef{#3}}}

\makenomenclature
\RequirePackage{ifthen}

\usepackage{longtable}
\usepackage{url}
\usepackage{floatrow}
\newcommand{\ones}{\mathbf{1}}
\newcommand{\mat}{\mathbf}

\newcommand\myeq[1]{\mathrel{\overset{\makebox[0pt]{#1}}{=}}}
\newcommand\mylt[1]{\mathrel{\overset{\makebox[0pt]{#1}}{\le}}}

\newcommand{\norm}[1]{\left\lVert#1\right\rVert}

\newtheorem{defi}{Definition}

\newtheorem{cor}{Corollary}

\begin{document}
\nomtypeD[a0]{\( \mat A\)}{Adjacency matrix }

\nomtypeD[a2]{\( d_i\)}{\(\displaystyle
\sum_j A_{ij} \), Degree of node $i.$ }

\nomtypeD[a3]{\( \mat d\)}{\(\displaystyle
\mat A \ones \), Vector of degrees }

\nomtypeD[a4]{\( \mat D\)}{\(\displaystyle
\text{diag}(\mat d) \), Diagonal matrix with degrees on the diagonal}

\nomtypeD[a6]{\( \mat P\)}{ $\mat A \mat D^{-1},$ Column-stochastic Markov matrix}

\nomtypeD[a7]{\( \mat v\)}{Preference vector of PageRank}

\nomtypeD[a8]{\( \widetilde{\mat P}\)}{\(\displaystyle
\alpha \mat P + (1-\alpha)\mat v \ones^T \) Transition Matrix of PageRank}

\nomtypeD[a9]{\( \mat Q\)}{\(\displaystyle
\mat D^{-1/2} \mat A \mat D^{-1/2} \), Symmetrized Markov matrix}

\nomtypeD[b1]{\( \overline{\mat Q}\)}{\(\displaystyle
\mat W^{-1/2} \overline{\mat A} \mat W^{-1/2} \), Average version of $\mat Q$}

\nomtypeD[a5]{\( \mat W\)}{\(\displaystyle
\mathbb E(\mat D) \), Expected Diagonal Matrix}

\nomtypeD[b2]{\( \mat R\)}{\(\displaystyle
(\mat I - \alpha \mat P)^{-1} \), Resolvent Matrix of $\mat P$}

\nomtypeD[b3]{\( \mat S\)}{\(\displaystyle
(\mat I - \alpha \mat Q)^{-1} \), Resolvent Matrix of $\mat Q$}

\nomtypeD[b5]{\( \mat u_1\)}{\(\displaystyle
{\mat D^{1/2}\ones}/{(\ones^T\mat D\mat 1)} \), Frobenius Perron eigenvector $\mat Q$}

\nomtypeD[b6]{\( \mat v^{'}\)}{\(\displaystyle
{n \mat D^{-1/2}\mat v} \), Degree-normalised and scaled $\mat v$}

\nomtypeD[b7]{\( \widetilde{\mat Q}\)}{\(\displaystyle
\mat Q  - \mat u_1 \mat u_1^T\), Projection of $\mat Q$ onto the orthogonal subspace of $\mat u_1$}
\nomtypeD[a1]{\( \overline{\mat A}\)}{\(\displaystyle
\mathbb E(\mat A)\), Expectation of adjacency matrix}

\nomtypeD[b4]{\( \widetilde{\mat R}\)}{\(\displaystyle
\mat W^{1/2} \left (  \mat{I} - \alpha \mat Q \right )^{-1} \mat W^{-1/2}\), Scaled resolvent matrix}


\mainmatter  

\title{PageRank in Undirected Random Graphs}

\titlerunning{PageRank in Undirected Random Graphs}

%
%

\author{K. Avrachenkov$^1$, A. Kadavankandy$^1$\thanks{Primary author, {\tt arun.kadavankandy@inria.fr}}, \\ L. Ostroumova Prokhorenkova$^{2,3}$ and A. Raigorodskii$^{2,3}$}

\authorrunning{K. Avrachenkov et al.}

\institute{Inria Sophia Antipolis, France \and
Yandex, Russia \and
Moscow Institute of Physics and Technology, Russia}

\maketitle

%
%

\begin{abstract}
PageRank has numerous applications in information retrieval, reputation systems, machine learning, and graph partitioning.
In this paper, we study PageRank in undirected random graphs with an expansion property. The Chung-Lu random graph is
an example of such a graph. We show that in the limit, as the size of the graph goes to infinity, PageRank can be approximated by a mixture of the restart distribution and the vertex degree distribution. We also extend the result to Stochastic Block Model (SBM) graphs, where we show that there is a correction term that depends on the community partitioning.
\keywords{PageRank, undirected random graphs, expander graphs,
Chung-Lu random graphs, Stochastic Block Model}
\end{abstract}

\vspace{-0.5cm}
\section{Introduction}
\vspace{-0.1cm}

PageRank has numerous applications in information retrieval \cite{H02,PBMW97,Yetal09}, reputation systems \cite{Getal13,KSG03},
machine learning \cite{Aetal08,AGMS12}, and graph partitioning
\cite{ACL06,C09}.
A large complex network can often be conveniently modeled by
a random graph. It is surprising that not many analytic studies are available for
PageRank in random graph models. We mention the work \cite{AL06} where PageRank
was analysed in preferential attachment models and the more recent works \cite{CLO14a,CLO14b},
where PageRank was analysed in directed configuration models. According to several studies \cite{Detal02,Fetal08,LSV07,VL10}, PageRank and in-degree
are strongly correlated in directed networks such as the Web graph.
Apart from some empirical
studies~\cite{B13,Perra}, to the best of our knowledge, there is no rigorous analysis of PageRank
on basic undirected random graph models such as the Erd\H{o}s-R\'enyi graph \cite{ER59} or the Chung-Lu graph \cite{CL02}.
In this paper, we attempt to fill this gap and show that under certain conditions on the preference vector and the spectrum of the graphs, PageRank in these models can be approximated by a mixture
of the preference vector and the vertex degree distribution when the size of the graph
goes to infinity. First, we show the convergence in total variation norm for a general family of random graphs with expansion property. Then, we specialize the results for the Chung-Lu
random graph model proving the element-wise convergence. We also analyse the asymptotics of PageRank on Stochastic Block Model (SBM) graphs, which are random graph models used to benchmark community detection algorithms \cite{holland1983}. In these graphs the asymptotic expression for PageRank contains an additional correction term that depends on the community partitioning. This demonstrates that PageRank captures properties of the graph not visible in the stationary distribution of a simple random walk.We conclude the paper with numerical experiments and several future research directions.

\section{Definitions}

Let $G^{(n)}=(V^{(n)},E^{(n)})$ denote a family of random graphs, where $V^{(n)}$ is a vertex set, $|V^{(n)}|=n$,
and $E^{(n)}$ is an edge set, $|E^{(n)}|=m$. Matrices and vectors related to the graph are denoted by bold letters, while their components are denoted by non-bold letters. We denote  by $\mat A^{(n)}$ the associated adjacency matrix
with elements
\[
A^{(n)}_{ij} =
\left\{ \begin{array}{ll}
1, & \mbox{if} \ i \ \mbox{and} \ j \ \mbox{are connected},\\
0, & \mbox{otherwise},
\end{array}\right.
\]
In the interest of compactness of notation, the superscript $n$ is dropped when it is not likely to cause confusion. In this work, since we analyze PageRank on undirected graphs, we have $\mat A^T=\mat A$. The personalized PageRank vector is denoted by $\boldsymbol{\pi}.$
We consider unweighted graphs; however our analysis easily extends to some families of weighted undirected graphs.
 Let $\ones$ be a column vector of $n$ ones and let $\mat d=\mat A\ones$ be the vector of degrees. It is helpful to define $\mat D=\mbox{diag}(\mat d)$,
a diagonal matrix with the degree sequence on its diagonal.

Let $\mat P=\mat A \mat D^{-1}$ be column-stochastic Markov transition matrix corresponding to the standard random walk on the graph
and let $\mat Q=\mat D^{-1/2}\mat A \mat D^{-1/2}$ be the symmetrized transition matrix, whose eigenvalues are the same as those of $\mat P.$
Note that the symmetrized transition matrix
is closely related to the normalized Laplacian ${\boldsymbol{\cal L}}=\mat I-\mat D^{-1/2}\mat A\mat D^{-1/2}=\mat I-\mat Q$ \cite{C97}, where $\mat I$ is the identity matrix. Further we will
also use the resolvent matrix $\mat R=[\mat I-\alpha \mat P]^{-1}$ and the symmetrized resolvent matrix $\mat S=[\mat I-\alpha \mat Q]^{-1}$.

Note that since $\mat Q$ is a symmetric matrix, its eigenvalues $\lambda_i,$ $i=1,...,n$ are
real and can be arranged in decreasing order,
i.e., $\lambda_1 \ge \lambda_2 \ge ...$ .
In particular, we have $\lambda_1=1$.
The value $\delta = 1-\max\{|\lambda_2|,|\lambda_n|\}$ is called the spectral gap.

In what follows, let $K,C$ be  arbitrary constants independent of graph size $n,$ which may change from one line to the next
(of course, not causing any inconsistencies).

For two functions $f(n),g(n),$ $g(n) = O(f(n))$ if $\exists C,N$ such that $\left| \frac{g(n)}{f(n)} \right| \leq C,$  $\forall n>N$ and $g(n) = {o}(f(n))$ if $\limsup_{n \to \infty} \left| \frac{g(n)}{f(n)} \right| = 0.$ Also $f(n) = \mathrm\omega(g(n))$ or $f(n) \gg g(n) $ if $g(n) = {o}(f(n)).$

We use $\mathbb P, \mathbb E$ to denote probability and expectation respectively. An event $E$ is said to hold with high probability (w.h.p.)
if $\exists N$ such that (s.t.) $\mathbb{P}(E) \geq 1 - {O}(n^{-c})$ for some $c > 0,$ $\forall n>N.$ Recall that if a finite number of events hold true w.h.p., then so does their intersection. Furthermore, we say that a sequence of random variables $X_n = {o}(1)$ w.h.p. if there exists a function $\psi(n) = {o}(1)$ such that the event $\{X_n \leq \psi(n)\}$ holds w.h.p.

In the first part of the paper, we study the asymptotics of PageRank for a family of random graphs with
the following two properties:

\smallskip
\begin{property}\label{prop:bounded_degrees}
For some $K$ w.h.p.,
$d^{(n)}_{max}/d^{(n)}_{min} \le K,$
where $d^{(n)}_{max}$ and $d^{(n)}_{min}$ are the maximum and minimum degrees, respectively.
\end{property}
\begin{property}\label{prop:fast_mixing} W.h.p.,
$\max\{|\lambda^{(n)}_2|,|\lambda^{(n)}_n|\}  = o(1).$
\end{property}

The above two properties can be regarded as a variation of the expansion property. In the standard case of an expander family,
one requires the graphs to be regular and the spectral gap
$\delta = 1-\max\{|\lambda_2|,|\lambda_n|\}$ to be bounded
away from zero (see, e.g., \cite{V12}). Property 1 is a relaxation of the regularity condition, whereas Property 2 is stronger than the requirement for the spectral gap to be bounded away from zero. These two properties allow us to consider several standard families of random graphs such as ER graphs, regular random graphs with increasing average degrees, and
Chung-Lu graphs. For Chung-Lu graphs Property 1 imposes some restriction on the degree spread of the graph.

\vspace{5pt}
\noindent \textit{Remark:} Property 2 implies that the graph is connected w.h.p., since the spectral gap is strictly greater than zero.

\vspace{5pt}

Later, we study the asymptotics of PageRank for specific classes of random graphs namely the Chung-Lu graphs, and the Stochastic Block Model.
 Recall that the Personalized PageRank vector with preference vector $\mat v$ is defined as the stationary
distribution of a modified Markov chain with transition matrix
\begin{equation}\label{eq:tildeP}
\widetilde{\mat P} = \alpha \mat P + (1-\alpha) \mat v \ones^T,
\end{equation}
where $\alpha$ is the so-called damping factor \cite{H02}. In other words, $\boldsymbol{\pi}$ satisfies
\begin{equation}
\boldsymbol{\pi} = \widetilde{\mat P} \boldsymbol{\pi} \label{eq:PRstat}, \\
\end{equation}
or,
\begin{equation}
\boldsymbol{\pi} = (1-\alpha) [\mat I-\alpha \mat P]^{-1} \mat v  = (1-\alpha) \mat R \mat v \label{eq:PPRexpl},
\end{equation}
where (\ref{eq:PPRexpl}) holds when $\alpha < 1.$
\vspace{-0.2cm}
\section{Convergence in total variation}\label{sec:tvconv}
\vspace{-0.1cm}

We recall that for two discrete probability distributions $u$ and $v$, the total variation distance $d_{\text{TV}}(u,v)$ is defined as $d_{\text{TV}}(u,v) = \frac{1}{2}\sum_i |u_i - v_i|.$ This can also be thought of as the $L^1$-norm distance measure in the space of probability vectors, wherein for $\mat x \in \mathbb{R}^{n},$ the $L^1$-norm is defined as $\norm{\mat x}_1 = \sum_i |x_i|.$ Since for any probability vector $\boldsymbol{\pi}, \text{ } \norm{\boldsymbol{\pi}}_1 = 1$ $\forall n,$ it makes sense to talk about convergence in 1-norm or TV-distance. Also recall that for a vector $\mat x \in \mathbb{R}^n,$ $\norm{\mat x}_2 = \sqrt{\sum_i |x_i|^2}$ is the $L^2$-norm. Now we are in a position to formulate our first result.

\begin{theorem}\label{th:1norm}
\label{total variation convergence}
Let a family of graphs $G^{(n)}$ satisfy Properties 1 and 2.
If, in addition, $\norm{\mat v}_2 = O(1/\sqrt{n})$, PageRank can be asymptotically approximated in total variation norm by a mixture of the restart distribution $\mat v$ and the vertex degree distribution. Namely, w.h.p.,
$$
d_{TV}(\boldsymbol{\pi}^{(n)},\overline{\boldsymbol{\pi}}^{(n)}) = o(1) \,\, \text{ as } n \to \infty,
$$
where
\begin{equation}
\label{eq:asymexpr}
\overline{\boldsymbol{\pi}}^{(n)}= \frac{\alpha \mat d^{(n)}}{\textnormal{vol}(G^{(n)})} +(1-\alpha) \mat v,
\end{equation}
with $\textnormal{vol}(G^{(n)}) = \sum_i d_i^{(n)}$.
\end{theorem}
\vspace{2pt}
\textit{Observations:}
\vspace{-2pt}
\begin{enumerate}
 \item This result says that PageRank vector asymptotically behaves like a convex combination of the preference vector and the stationary vector of a standard random walk with transition matrix $\mat P;$ with the weight being $\alpha,$ and that it starts to resemble the random walk stationary vector as $\alpha$ gets close to $1.$
\item One of the possible intuitive explanations of the result of Theorem \ref{th:1norm} is based on the observation that when Properties 1 \& 2 hold, as $n \to \infty,$ the random walk mixes approximately in one step and so for any probability vector $\mat x$ $\mat{P} \mat x$ is roughly equal to $\mat d/\mbox{vol}(G),$ the stationary distribution of the simple random walk. The proposed asymptotic approximation for PageRank can then be seen to follow from the series representation of PageRank if we replace $\mat P \mat v$ by $\mat d/\mbox{vol}(G).$ Note that since $\mat d/\textnormal{vol}(G)$ is the stationary vector of the simple random walk, if $\mat P \mat v = \mat d/\textnormal{vol}(G),$ it also holds that $\mat P^k \mat v = \mat d/\textnormal{vol}(G), \forall k \ge 2.$  Making these substitutions in the series representation of PageRank, namely
\begin{equation}
\label{eq:forpi}
\boldsymbol{\pi} = (1 - \alpha) \left( \mat I + \alpha \mat P + \alpha^2 \mat P^2 + \ldots   \right) \mat v,
\end{equation}
we obtain
\begin{align*}
\boldsymbol{\pi} &= (1-\alpha) \mat v + (1-\alpha)\alpha (1 + \alpha + \alpha^2 + \ldots)\frac{\mat d}{\textnormal{vol}(G)}\\
&= (1-\alpha) \mat v + \alpha \frac{\mat d}{\textnormal{vol}(G)}.
\end{align*}
\item The condition on the 2-norm of the preference vector $\mat v$ can be viewed as a constraint on its allowed localization.
\end{enumerate}

\textbf{Proof of Theorem \ref{th:1norm}:}
First observe from (\ref{eq:tildeP}) that when $\alpha = 0,$ we have $\widetilde{\mat{P}} = \mat v \ones^T,$ hence from (\ref{eq:PRstat}) we obtain $\boldsymbol{\pi} = \mat v,$ since $\ones^T \boldsymbol{\pi} = 1.$ Similarly for the case $\alpha = 1,$ $\widetilde{\mat{P}} = \mat P$ and so $\boldsymbol{\pi}$ in this case is just the stationary distribution of the original random walk, which is well-defined and equals $\frac{\mat d}{\textnormal{vol}(G)}$ since by Property 2 the graph is connected. Examining (\ref{eq:asymexpr}) for these two cases we can see that the statement of the theorem holds trivially for both $\alpha = 0$ and $\alpha = 1.$ In what follows, we consider the case $0 < \alpha < 1.$
We first note that the matrix $\mat Q = \mat D^{-1/2}\mat A \mat D^{-1/2}$ can be written as follows by Spectral Decomposition Theorem \cite{BhatiaSpr}:
\begin{equation}
\label{eq:Qeig}
\mat Q = \mat u_1 \mat u_1^T + \sum_{i=2}^n \lambda_i \mat u_i \mat u_i^T,
\end{equation}
where $1 = \lambda_1 \geq \lambda_2 \geq \ldots \geq \lambda_n$ are the eigenvalues and $\{\mat u_1,\mat u_2,  \ldots \mat u_n \}$ with $\mat u_i \in \mathbf{R}^n$ and $\norm{\mat u_i}_2 = 1$ are the corresponding orthogonal eigenvectors of $\mat Q.$ Recall that $\mat u_1=\mat D^{1/2}\ones/\sqrt{\ones^T \mat D\ones}$ is the Perron--Frobenius eigenvector.
Next, we rewrite (\ref{eq:PPRexpl}) in terms of the matrix $\mat Q$ as follows
\begin{equation}\label{eq:expprank}
\boldsymbol{\pi} = (1-\alpha) \mat D^{1/2} [\mat I-\alpha \mat Q]^{-1} \mat D^{-1/2} \mat v.
\end{equation}
Substituting (\ref{eq:Qeig}) into (\ref{eq:expprank}), we obtain
\begin{align*}
\boldsymbol{\pi} &
= (1-\alpha) \mat D^{1/2}
\left(\frac{1}{1-\alpha}\mat u_1 \mat u_1^T + \sum_{i=2}^n \frac{1}{1-\alpha\lambda_i} \mat u_i \mat u_i^T \right) \mat D^{-1/2} \mat v  \\
& = \mat D^{1/2}\mat u_1 \mat u_1^T \mat D^{-1/2}\mat v + (1 - \alpha) \mat D^{1/2} \left( \sum_{i\neq 1} \frac{1}{1 - \alpha \lambda_i} \mat u_i \mat u_i^T \right)\mat D^{-1/2} \mat v. \\
\end{align*}

Let us denote the error vector by $\boldsymbol{\epsilon} = \boldsymbol{\pi} - \overline{\boldsymbol{\pi}}$. Note that since $\mat u_1 = \frac{\mat D^{1/2} \ones}{\sqrt{\textnormal{vol}(G)}},$ we can write $\overline{\boldsymbol{\pi}}$ as
\begin{align*}
\overline{\boldsymbol{\pi}} &= \alpha \frac{\mat d}{\textnormal{vol}(G)} + (1-\alpha) \mat v \\
&\myeq{(a)} \alpha \frac{\mat D \ones \ones^T \mat v}{\textnormal{vol}(G)} + (1-\alpha) \mat D^{1/2}\mat D^{-1/2} \mat v\\
&= \alpha \mat D^{1/2} \frac{\mat D^{1/2} \ones}{\sqrt{\textnormal{vol}(G)}} \frac{\ones^T\mat D^{1/2}}{\sqrt{\textnormal{vol}(G)}} \mat D^{-1/2} \mat v + (1-\alpha) \mat D^{1/2}\mat D^{-1/2} \mat v\\
&= \alpha \mat D^{1/2}\mat u_1 \mat u_1^T \mat D^{-1/2}\mat v + (1-\alpha) \mat D^{1/2}\mat D^{-1/2}\mat v,
\end{align*}
where in (a) above we used the fact that $\ones^T \mat v = 1,$ since $\mat v$ is a probability vector.
Then, we can write $\boldsymbol{\epsilon}$ as
\begin{align}
\boldsymbol{\epsilon}
&= \boldsymbol{\pi} - \alpha \mat D^{1/2}\mat u_1 \mat u_1^T \mat D^{-1/2} \mat v -  (1 - \alpha) \mat D^{1/2} \mat I \mat D^{-1/2}\mat v \nonumber \\
 &= (1 - \alpha) \mat D^{1/2} \left ( \sum_{i \neq 1} \frac{\mat u_i \mat u_i^T}{1 - \alpha \lambda_i} - (\mat I - \mat u_1 \mat u_1^T)\right )  \mat D^{-1/2} \mat v \nonumber \\
&= (1 - \alpha) \mat D^{1/2} \left ( \sum_{i \neq 1} \mat u_i \mat u_i^T \frac{\alpha \lambda_i}{1 - \alpha \lambda_i} \right)  \mat D^{-1/2} \mat v. \label{eq:epsilonerr}
\end{align}
Now let us bound the $L^1$-norm $\norm{\boldsymbol{\epsilon}}_1$ of the error:
\begin{align}
\norm{\boldsymbol{\epsilon} }_1 /(1 - \alpha)
& \mylt{(a)} \sqrt{n} \|\boldsymbol{\epsilon}\|_2/(1-\alpha) \nonumber \\
& \mylt{(b)} \sqrt{n}\|\mat D^{1/2}\|_2 \left\| \sum_{i \neq 1} \mat u_i \mat u_i^T \frac{\alpha \lambda_i}{1 - \alpha \lambda_i} \right\|_2\|\mat D^{-1/2}\|_2 \|\mat v\|_2 \nonumber \\
&\mylt{(c)} \sqrt{d_{max}/d_{min}} \sqrt{n} \max_{i > 1}\left| \frac{\alpha \lambda_i}{1 - \alpha \lambda_i}\right| \norm{\mat v}_2 \nonumber \\
& \le C  \sqrt{d_{max}/d_{min}}  \max(|\lambda_2|, |\lambda_{n}|) \label{normbound}
\end{align}
where in (a) we used the fact that for any vector $\mat{x} \in \mathbb R^n,$ $\|\mat x\|_1 \le \sqrt{n} \|\mat x\|_2$ by Cauchy-Schwartz inequality. In (b) we used the submultiplicative property of matrix norms, i.e., $\norm{\mat A \mat B}_2 \leq \norm{\mat A}_2 \norm{\mat B}_2$. We obtain (c) by noting that the norm of a diagonal matrix is the leading diagonal value and the fact that for a symmetric matrix the 2-norm is the largest eigenvalue in magnitude. The last inequality is obtained by noting that the assumption $\lambda_i = o(1) $ w.h.p. $\forall i>1$ implies that $\exists N$ s.t. $\forall n > N,$ $|1-\alpha\lambda_i| > C$ for some constant C and the fact that $\norm{\mat v}_2 = O(1/\sqrt{n}).$

\noindent Observing that $d_{max}/d_{min}$ is bounded w.h.p. by Property~1 and $\max(|\lambda_2|, |\lambda_{n}|) = o(1)$ w.h.p. by Property~2 we obtain our result.  \qed
\smallskip

Note that in the case of standard PageRank, $v_i=1/n, 1\le i\le n,$ and
hence $\norm{\mat v}_2 = O(1/\sqrt{n}),$ but Theorem~\ref{th:1norm} also admits more general preference vectors than the uniform one.

\begin{cor}
The statement of Theorem~\ref{th:1norm} also holds with respect to the weak convergence, i.e., for any function $f$ on $V$ such that $\max_{x \in V} |f(x)| \leq 1,$ $$ \sup \left \{ \sum_v f(v) \pi_v - \sum_v f(v) \overline{\pi}_v \right \} = o(1) \quad \mbox{w.h.p.} $$
\end{cor}
\textbf{Proof:}
This follows from Theorem~\ref{th:1norm} and the fact that the left-hand side of the above equation is upper bounded by
$2\, d_{\text{TV}}(\boldsymbol{\pi}_n , \overline{\boldsymbol{\pi}}_n)$ \cite{AsherPeres2009}.
\qed
\vspace{-0.2cm}
\section{Chung-Lu random graphs}
\vspace{-0.1cm}

In this section, we study the PageRank for the Chung-Lu model~\cite{CL02} of random graphs. These results naturally hold for ER graphs also. The spectral properties of Chung-Lu graphs have been studied extensively in a series of papers by Fan Chung et al \cite{CLV03,Chung2011}.

\subsection{Chung-Lu Random Graph Model}

Let us first provide a definition of the Chung-Lu random graph model.

\begin{defi}{\bf Chung-Lu Random Graph Model}
A Chung-Lu graph $\mathcal{G}(w)$ with an expected degree vector $\mat w = (w_1,w_2, \ldots w_n)$, where $w_i$ are positive real numbers, is generated by drawing an edge between any two vertices $v_i$ and $v_j$ independently of all other pairs, with probability $p_{ij} = \frac{w_i w_j}{\sum_k w_k}.$ To ensure that the probabilities $p_{ij}$ are well-defined, we need  $\max_i w_i^2 \leq \sum_k w_k$.
\end{defi}
In the following, let $w_{\max} = \max_i w_i$ and $w_{\min} = \min_i w_i.$
Below we specify a corollary of Theorem \ref{total variation convergence} as applied to these graphs. But before that we need the following lemmas about Chung-Lu graphs mainly taken from \cite{CLV03,Chung2011}.

\begin{lemma}
\label{degree_concentration}
If the expected degrees $w_1,w_2,\ldots w_n$ satisfy $w_{\min} \gg \log(n),$ then in $\mathcal{G}(w)$ we have, w.h.p.,
$\max_i |\frac{d_i}{w_i} - 1| = o(1)$.
\end{lemma}

In the proof we use Bernstein Concentration Lemma~\cite{Bilsley}:
\begin{lemma}\label{lemma:berny}(Bernstein Concentration Lemma~\cite{Bilsley})
If $Y_n = X_1 + X_2 + \ldots X_n,$ where $X_i$ are independent random variables such that $|X_i|\le b$ and if $B_n^2 = \mathbb{E}(Y_n - \mathbb{E}(Y_n))^2 $ then
\[
\mathbb{P}\{|Y_n - \mathbb{E}(Y_n)| \geq \epsilon\} \leq 2 \exp \frac{-\epsilon^2}{2(B_n^2 + b\epsilon/3 )},
\]
for any $\epsilon > 0.$
\end{lemma}
\textbf{Proof of Lemma \ref{degree_concentration}:}
This result is shown in the sense of convergence in probability in the proof of \cite[Theorem~2]{Chung2011}; using Lemma \ref{lemma:berny} we show the result holds w.h.p.
By a straight forward application of Lemma \ref{lemma:berny} to the degrees $d_i$ of the Chung-Lu graph we obtain
$$
\mathbb{P}\left(\max_{1 \leq i \leq n} \left|\frac{d_i}{w_i} - 1 \right| \geq \beta \right) \leq \frac{2}{n^{c/4 - 1}}, \quad \mbox{if} \quad \beta \geq \sqrt{\frac{c\log(n)}{w_{\min}}}=o(1)
$$
if $w_{\min} \gg \log(n)$.\qed
We present below a perturbation result for the eigenvalues of Hermitian matrices, called Weyl's inequalities, which we will need for our proofs.
\begin{lemma} \cite[Theorem ~4.3.1]{hornjohn12}
\label{lem:weyl}
Let $\mat A,\mat B \in \mathbb R^{n\times n}$ be Hermitian and let the eigenvalues $\lambda_i(\mat A),$ $\lambda_i(\mat B)$ and $\lambda_i(\mat A+\mat B)$ be arranged in decreasing order. For each $k = 1,2, \ldots n$ we have
\[
|\lambda_k(\mat A+\mat B) - \lambda_k(\mat A)| \le \|\mat B\|_2,
\]
where $\|\mat B\|_2$ is the induced 2-norm or the spectral norm of $\mat B.$
\end{lemma}
The following lemma is an application of Theorem 5 in \cite{CLV03}.

\begin{lemma}
\label{twonormchunglu}
If $w_{\max} \leq K w_{\min},$ for some $K>0$ and $\overline{w} = \sum_k w_k/n \gg \log^6(n)$, then for $\mathcal{G}(w)$ we have almost surely (a.s.) $$\norm{\mat C}_2 = \frac{2}{\sqrt{\overline{w}}}(1 + o(1)),$$ where $ \mat C = \mat W^{-1/2}\mat A \mat W^{-1/2} - \boldsymbol{\chi}^{T} \boldsymbol{\chi}$, $\mat W = \mbox{diag}(\mat{w}),$ and $\mat \chi_i = \sqrt{w_i/\sum_k w_k}$ is a row vector.
\end{lemma}
\textbf{Proof:}
 It can be verified that when $w_{\max} \leq K w_{\min}$ and $\overline{w} \gg \log^6(n),$ the condition in \cite[Theorem~5]{CLV03}, namely, $w_{\min} \gg \sqrt{\overline{w}} \log^3(n),$ is satisfied and hence the result follows.\qed

\begin{lemma}
\label{second eigen chung}
For $\mathcal{G}(w)$ with $w_{\max}\leq K w_{\min},$ and  $\overline{w}\gg \log^6(n),$ $$\max(\lambda_2(\mat P),-\lambda_{n}(\mat P)) = o(1) \quad \mbox{w.h.p.},$$  where $\mat P$ is Markov matrix.
\end{lemma}
{\bf Proof:}
Recall that $\mat Q = \mat D^{-1/2} \mat A \mat D^{-1/2}$ is the normalized adjacency matrix. We want to be able to bound the eigenvalues $\lambda_i, i \ge 2$ of $\mat Q.$ We do this in two steps. Using Lemmas \ref{degree_concentration} and \ref{lem:weyl} we first show that if we replace the degree matrix $\mat D$ in the expression for $\mat Q$ by the expected degree matrix $\mat W = \mathbb{E}(\mat D),$ the eigenvalues of the resulting matrix are close to those of $\mat Q.$ Then, using Lemma \ref{twonormchunglu} we show that the eigenvalues of $\mat W^{-1/2}\mat A \mat W^{-1/2}$ roughly coincide with those of $\boldsymbol{\chi}^{T} \boldsymbol{\chi},$ which is a unit rank matrix and hence only has a single non-zero eigenvalue. Thus we arrive at the result of Lemma \ref{second eigen chung}. Now we give the detailed proof.

The first step, $\|\mat Q - \mat W^{-1/2}\mat A \mat W^{-1/2}\|_2 = o(1)$ w.h.p. follows from Lemma \ref{degree_concentration} and the same argument as in the last part of the proof of Theorem 2 in \cite{Chung2011}. We present the steps in the derivation here for the sake of completeness.

Since the 2-norm of a diagonal matrix is the maximum diagonal in absolute value, we have

\begin{equation}\label{eq:normerr}
 \|\mat W^{-1/2}\mat D^{1/2} - \mat I \|_2 = \max_{\{i = 1,2,\ldots \}} \left|\sqrt{\frac{d_i}{w_i}} - 1\right| \le \max_{\{i = 1,2,\ldots \}} \left|{\frac{d_i}{w_i}} - 1 \right|  = o(1),
\end{equation}
by Lemma \ref{degree_concentration}. Also observe that
\begin{equation}
\label{eqnormq}
\|\mat Q\|_2 = \max_{\{i=1,2,\ldots n \}} |\lambda_i(\mat Q)| = \max_{\{i=1,2,\ldots n \}} |\lambda_i(\mat P)| = 1.
\end{equation}
We now proceed to bound the norm of the difference $\|\mat Q - \mat W^{-1/2}\mat A \mat W^{-1/2}\|$ as follows
\begin{eqnarray}
\lefteqn{\|\mat Q - \mat W^{-1/2}\mat A \mat W^{-1/2}\|_2} \nonumber\\
&=& \| \mat Q -  \mat W^{-1/2} \mat D^{1/2} \mat D^{-1/2}\mat A \mat D^{-1/2} \mat D^{1/2} \mat W^{-1/2}\|_2 \nonumber\\
&=& \|\mat Q - \mat W^{-1/2} \mat D^{1/2} \mat Q  \mat D^{1/2} \mat W^{-1/2} \|_2 \nonumber\\
&=& \|\mat Q - \mat W^{-1/2} \mat D^{1/2} \mat Q  + \mat W^{-1/2} \mat D^{1/2} \mat Q -  \mat W^{-1/2} \mat D^{1/2} \mat Q \mat D^{1/2} \mat W^{-1/2}\|_2 \nonumber \\
&\myeq{(a)}& \|(\mat I  - \mat W^{-1/2} \mat D^{1/2}) \mat Q\|_2 + \|\mat W^{-1/2} \mat D^{1/2} \mat Q(\mat I - \mat D^{1/2} \mat W^{-1/2})\|_2 \nonumber \\
&\mylt{(b)}& \|(\mat I  - \mat W^{-1/2} \mat D^{1/2})\|_2 \|\mat Q\|_2 + \|\mat W^{-1/2} \mat D^{1/2}\|_2 \|\mat Q\|_2 \|\mat I - \mat D^{1/2} \mat W^{-1/2}\|_2 \nonumber \\
&\myeq{(c)}& o(1) + (1 + o(1))o(1) = o(1) \quad w.h.p., \label{eq:boundonq}
\end{eqnarray}
where (a) follows from triangular inequality of norms, in (b) we used submultiplicativity of matrix norms, and (c) follows from (\ref{eq:normerr}), (\ref{eqnormq}) and the fact that  $\|\mat W^{-1/2} \mat D^{1/2}\|_2 \le \|\mat I \|_2 + \|\mat W^{-1/2} \mat D^{1/2} - \mat I \|_2 = (1 + o(1)).$

By Lemma \ref{lem:weyl} we have for any $i,$
\begin{equation}\label{eq:fpart}
|\lambda_i(\mat Q) - \lambda_i(\mat W^{-1/2} \mat A \mat W^{-1/2})| \le \|\mat Q - \mat W^{-1/2} \mat A \mat W^{-1/2}\|_2 = o(1),
\end{equation}
by (\ref{eq:boundonq}). Furthermore, using Lemma \ref{lem:weyl} and the fact that $\lambda_i(\boldsymbol{\chi}^T \boldsymbol{\chi}) = 0$ for $i>1,$ we have for $i \ge 2,$
\begin{eqnarray}
\lefteqn{|\lambda_i(\mat W^{-1/2} \mat A \mat W^{-1/2})|}\nonumber \\
&=& |\lambda_i(\mat W^{-1/2} \mat A \mat W^{-1/2}) - \lambda_i(\boldsymbol{\chi}^T \boldsymbol{\chi})|
\le \|\mat W^{-1/2} \mat A \mat W^{-1/2} - \boldsymbol{\chi}^T \boldsymbol{\chi}\|_2 \nonumber \\
 &=& o(1),  \label{eq:diffbnd}
\end{eqnarray}
where the last inequality follows from Lemma \ref{twonormchunglu}.\\
Now recall that $\max(\lambda_2(\mat P),-\lambda_{n}(\mat P)) = \max_{\{i \ge 2\}}|\lambda_i(\mat Q)|.$
We have for any $i,$
\begin{align}
|\lambda_i(\mat Q)| \le |\lambda_i(\mat Q) - \lambda_i(\mat W^{-1/2} \mat A \mat W^{-1/2})| + ||\lambda_i(\mat W^{-1/2} \mat A \mat W^{-1/2})|,
\end{align}
which implies from (\ref{eq:fpart}) and (\ref{eq:diffbnd}):
\[
\max_{\{i \ge 2\}} |\lambda_i(\mat Q)| = o(1).
\]
\qed
Armed with these lemmas we now present the following corollary of Theorem \ref{th:1norm} in the case of Chung-Lu graphs.
\begin{cor}
\label{prChungLuTV}
Let $\norm{\mat v}_2 = O(1/\sqrt{n}),$ and $\alpha \in (0,1).$ Then PageRank $\boldsymbol{\pi}$ of the Chung-Lu graph $\mathcal{G}(w)$ can asymptotically be approximated in TV distance by $\overline{\boldsymbol{\pi}},$ defined in Theorem \ref{total variation convergence}, if $\overline{w} \gg \log^6(n)$ and $w_{\max} \leq K w_{\min}$ for some $K$ that does not depend on $n.$
\end{cor}
\textbf{Proof:}
Using Lemma \ref{degree_concentration} and the condition that $w_{\max} \leq K w_{\min},$ one can show that $\exists K^{'}$ s.t. $\frac{d_{max}}{d_{min}} \leq K^{'}$ w.h.p. Then the result is a direct consequence of Lemma \ref{second eigen chung} and the inequality from~\eqref{normbound}.\qed

We further note that this result also holds for ER graphs $\mathcal{G}(n,p_n)$ with $n$ nodes and edge probability $p_n$ such that $n p_n \gg \log^6(n),$ where we have $(w_1,w_2,\ldots w_n) = (np_n, np_n, \ldots np_n).$

\subsection{Element-wise Convergence}

In Corollary \ref{prChungLuTV} we proved the convergence of PageRank in TV distance for Chung-Lu random graphs. Note that since each component of PageRank could decay to zero as the graph size grows to infinity, this does not necessarily guarantee convergence in an element-wise sense. In this section, we provide a proof for our convergence conjecture to include the element-wise convergence of the PageRank vector. Here we deviate slightly from the spectral decomposition technique and eigenvalue bounds used hitherto, and instead rely on well-known concentration bounds to bound the error in convergence.

Let $\overline{\boldsymbol{\Pi}} = \mbox{diag} \{ \overline{\pi}_1,\overline{\pi}_2, \ldots \overline{\pi}_n \}$ be a diagonal matrix whose diagonal elements are made of the components of
the approximated PageRank vector  and $\widetilde{\boldsymbol{\delta}} = \overline{\boldsymbol{\Pi}}^{-1}(\boldsymbol{\pi} - \overline{\boldsymbol{\pi}}),$ i.e., $\widetilde{\delta}_i = (\pi _i - \overline{\pi} _i)/\overline{\pi}_i = \epsilon_i/\overline{\pi}_i,$ where $\boldsymbol{\epsilon}$ is the unnormalized error defined in Section \ref{sec:tvconv}. Then using (\ref{eq:epsilonerr}) we obtain
\[
\widetilde{\delta}_i = \left((1 - \alpha)v_i + \alpha \frac{d_i}{\text{vol(G)}}\right)^{-1}  \left [\mat D^{1/2} \sum_{j \neq 1} \frac{\alpha \lambda_j}{1 - \alpha \lambda_j} \mat u_j \mat u_j^T \mat D^{-1/2} \mat v \right ]_i.
\]
Therefore, using $\mat v'$ to denote $n \mat D^{-1/2} \mat v$ we can bound $\norm{\widetilde{\boldsymbol{\delta}}}_{\infty} = \max_i |\widetilde{\delta}_i|$ as follows
\begin{align}
\norm{\widetilde{\boldsymbol{\delta}}}_{\infty} &\le \frac{1}{\min_i \left((1-\alpha)v_i + \alpha \frac{d_i}{\textnormal{vol}(G)}\right)}\left\|\mat D^{1/2} \sum_{j \neq 1} \frac{\alpha \lambda_j}{1 - \alpha \lambda_j} \mat u_j \mat u_j^T \mat D^{-1/2} \mat v\right\|_{\infty}\\
&\leq \frac{\sum_i d_i/n}{\alpha d_{\min}} \sqrt{d_{max}} \norm{\sum_{j \neq 1} \frac{\alpha \lambda_j}{1 - \alpha \lambda_j} \mat u_j \mat u_j^{T} \mat v'}_{\infty}\label{infty_normbound}.
\end{align}
Here $d_{\min}$ denotes $\min_i d_i.$ To obtain (\ref{infty_normbound}) we used the submultiplicativity property of matrix norms, the fact that $\|\mat D^{1/2}\|_{\infty} = \sqrt{\max_i d_i} = \sqrt{d_{\max}}$ and the fact that $v_i \ge 0, \forall i \in V.$

Define $\widetilde{\mat Q} = \mat Q - \mat u_1 \mat u_1^T,$ the restriction of the matrix $\mat Q$ to the orthogonal subspace of $\mat u_1.$

\begin{lemma}
\label{Bound on Sv}
For a Chung-Lu random graph $\mathcal{G}(w)$ with expected degrees $w_1,\ldots w_n$, where $w_{\max} \leq K w_{\min}$ and $w_{\min} \gg \log(n),$ we have w.h.p.,
\[
\norm{\widetilde{\mat Q}\mat v'}_{\infty} = o(1/\sqrt{w_{\min}}),
\]
when $v_i = O(1/n) \ \forall i.$
\end{lemma}
This lemma can be proven by a few applications of Lemma \ref{degree_concentration} and Bernstein's concentration inequality. To keep the train of thought intact, please refer to Appendix \ref{ap:lemmaSv} for a detailed proof of this lemma.

In the next lemma we prove an upper bound on the infinity norm of the matrix $\mat S = (\mat I - \alpha \mat Q)^{-1}.$

\begin{lemma}
\label{bound_inv_infty}
Under the conditions of Lemma \ref{Bound on Sv}, $\norm{\mat S}_{\infty} \leq C $  w.h.p., where $\mat C$ is a number independent of $n$ that depends only on $\alpha$ and $K$.
\end{lemma}
{\bf Proof:}
Note that $\mat S = (\mat I - \alpha \mat Q)^{-1} = \mat D^{-1/2}(\mat I - \alpha \mat P)^{-1} \mat D^{1/2}.$ Therefore, $\norm{\mat S}_{\infty} \leq \sqrt{\frac{d_{max}}{d_{min}}} \norm{(\mat I - \alpha \mat P)^{-1}}_{\infty}$ and the result follows since $\norm{(\mat I - \alpha \mat P)^{-1}}_{\infty} \leq \frac{1}{1 - \alpha}$ \cite{LangMeyer} and using Lemma~\ref{degree_concentration}.
\qed
Now we are in a position to present our main result in this section.

\begin{theorem}
\label{elementwise}
Let $v_i = O(1/n) \,\, \forall i,$ and $\alpha < 1.$ PageRank $\boldsymbol{\pi}$ converges element-wise to $\overline{\boldsymbol{\pi}} = (1 - \alpha) \mat v + \alpha \mat d/\textnormal{vol}(G),$ in the sense that $\max_i \ (\pi _i - \overline{\pi} _i)/\overline{\pi}_i = o(1)$ w.h.p., on the Chung-Lu graph $\mathcal G(w)$ with expected degrees $\{w_1,w_2,\ldots w_{n} \}$ such that $w_{\min} > \log^c(n)$ for some $c > 1$ and $w_{\max} \leq K w_{\min},$ for some $K,$ a constant independent of $n.$
\end{theorem}
\textbf{Proof:} Define $\mat Z = \sum_{i \neq 1} \frac{\alpha \lambda_i}{1 - \alpha \lambda_i} \mat u_i \mat u_i^T. $
We then have:

\begin{align}
 \mat Z &= \sum_{i = 1}^{n} \frac{\alpha \lambda_i}{1 - \alpha \lambda_i} \mat u_i \mat u_i^T - \frac{\alpha}{1 - \alpha} \mat u_1 \mat u_1^T \nonumber \\
 & = (\mat I - \alpha \mat Q)^{-1} \alpha \mat Q - \frac{\alpha}{1 - \alpha} \mat u_1 \mat u_1^T \nonumber \\
 & = \mat S \left [ \alpha \mat Q - \frac{\alpha}{1 - \alpha} (\mat I - \alpha \mat Q) \mat u_1 \mat u_1^T\right ] \nonumber \\
 & = \alpha \mat S \widetilde{\mat Q} \label{eq:z}
\end{align}

Now from (\ref{infty_normbound}) we have
\begin{align*}
\norm{\widetilde{\boldsymbol\delta}}_{\infty} &\le   C \frac{\sum_i d_i/n}{ d_{\min}} \sqrt{d_{max}} \|\mat S \widetilde{\mat Q}\mat v^{'}\|_{\infty} \\
& \mylt{(a)} C \frac{\sum_i d_i/n}{ d_{\min}} \sqrt{d_{max}}o(1/\sqrt{w_{\min}}) \nonumber \\
& \le C \frac{d_{\max}}{d_{\min}} \sqrt{w_{\max}(1+o(1))} \frac{1}{\sqrt{w_{\min}}}o(1)  \label{replacedw} \\
&= C \frac{w_{\max}}{w_{\min}} \sqrt{\frac{w_{\max}}{w_{\min}}} (1+o(1)) o(1) \\
&= C \left(\frac{w_{\max}}{w_{\min}}\right)^{\frac{3}{2}}o(1) \\
& \le C o(1) \quad  \mbox{w.h.p.},
\end{align*}
where in (a) we used (\ref{eq:z}) and Lemmas \ref{Bound on Sv} and \ref{bound_inv_infty}. The rest of the inequalities are obtained by repeatedly using the fact that $d_{\max} = w_{\max}(1+o(1))$ and $d_{min} = w_{\min}(1+o(1)),$ from Lemma \ref{degree_concentration}. The last step follows from the assumption that $w_{\max}\le Kw_{\min}$ for some constant $K.$\qed
\begin{corollary} [ER Graphs]\\
For an ER graph $\mathcal{G}(n,p_n)$ such that $np_n \gg \log(n),$ we have that asymptotically the personalized PageRank $\boldsymbol\pi$ converges pointwise to $\overline{\boldsymbol\pi}$ for $\mat v$ such that $v_i = O(1/n).$
\end{corollary}



\section{Asymptotic PageRank for the Stochastic Block Model}\label{sec:sbmprank}
In this section, we extend the analysis of PageRank to Stochastic Block Models (SBM) with constraints on average degrees.
The SBM is a random graph model that reflects the community structure prevalent in many online social networks. It was first introduced in \cite{holland1983} and has been analyzed subsequently in several works, specifically in the community detection literature, including \cite{condon2001},\cite{Karrer2011}, \cite{rohe2011}, \cite{ACK15} and several extensions thereof as in \cite{Heimlicher2012} and \cite{zhao2012}, and the references therein.

For the sake of simplicity we focus on an SBM graph with two communities, but the idea of the proof extends easily to generalizations of this simple model.
\begin{definition} \label{def:sbm}
[Stochastic Block Model (SBM) with two communities]: An SBM graph $\mathcal{G}(m,n-m,p,q)$ with two communities is an undirected graph on a set of disjoint vertices $C_1, C_2$ such that $ C_1 \cup C_2 = V,$ and let $|C_1| = m$ and $|C_2| = n-m$. Furthermore, if two vertices $i,j \in C_k, k = 1,2$, then $\mathbb{P}( (i,j) \in E ) = p$, if $i \in C_1$ and $j \in C_{2}$,  then $\mathbb{P}( (i,j) \in E ) = q.$ The probabilities $p,q$ may scale with $n$ and we assume that $m > n/2$ and $p > q;$ this last assumption is necessary for modeling the community structure of a network.
\end{definition}

\noindent \textit{Remark:} For the sake of simplicity, we assume that the edge probabilities within both communities are equal to $p,$ but this is a minor assumption and can be generalised so that community 1 has a different edge probability to community 2.

For an SBM graph we use $w_{\max}$ and $w_{\min}$ to denote the maximum and the minimum expected degrees of the nodes respectively. From Definition \ref{def:sbm}, by our assumption on $m,p$ and $q,$ we have $w_{\max} = m p+ (n-m) q$ and $w_{\min} = (n-m)p + m q.$ Note that our results only depend on these two parameters.
We present our main result on SBM graphs in the following theorem.

\begin{theorem}
\label{theo:sbm}
For a Stochastic Block Model with $w_{\min} = \mathrm{\omega}(\log^3(n))$ and $\frac{w_{\max}}{w_{\min}} \le C,$ PageRank with preference vector $\mat{v}$ such that $\|\mat{v}\|_2 = O(\frac{1}{\sqrt{n}})$ satisfies
\[
\| \boldsymbol{\pi} - \overline{\boldsymbol{\pi}}_{\textnormal{SBM}} \|_{\text{TV}} = o(1)
\]
w.h.p., where
\begin{equation}
\label{eq:prSBM}
\overline{\boldsymbol{\pi}}_{\textnormal{SBM}}  = (1- \alpha) \left( \mat I - \alpha \overline{\mat P} \right)^{-1} \mat{v}.
\end{equation}
Here $\overline{\mat P}$ represents the ``average'' Markov matrix given as $\overline{\mat P} = \overline{\mat A}  \mat W^{-1}$ where $\mat W= \mathbb{E}(\mat D)$  and $\overline{\mat A}= \mathbb{E}(\mat{A}).$
\end{theorem}

\noindent \textit{Discussion:} Let us look at the permissible values of $m,p,q$ under the assumptions in the above theorem. Recall that we have $w_{\min} = (n-m)p + m q > nq.$ Therefore the condition on the growth of minimum expected degree is met, for example, if $q = \omega(\log^3(n)/n).$ On the other hand we have
\[
\frac{w_{\max}}{w_{\min}} = \frac{mp + (n-m)q}{(n-m)p + m q} = \frac{\frac{m}{n-m} \frac{p}{q} + 1}{\frac{m}{n-m}+\frac{p}{q}} \quad,
\]
which remains bounded if either $m/(n-m)$ or $p/q$ tends to infinity, but not both.

The following corollary of Theorem \ref{theo:sbm} gives an interesting expression for PageRank for an SBM graph with two equal-sized communities.
\begin{corollary}
For an SBM graph as in Definition \ref{def:sbm}, with $m = n/2,$ (n assumed to be even) such that $p+q \gg \log^3(n)/n$ the PageRank vector $\boldsymbol{\pi}$ with preference vector $\mathbf{v}$ such that $\|\mathbf{v}\|_2 = O(\frac{1}{\sqrt{n}})$ satisfies
\[
\| \boldsymbol{\pi} - \overline{\boldsymbol{\pi}}_{\textnormal{SBM}} \|_{TV} \to 0
\]
w.h.p as $n \to \infty$ where
\begin{equation}\label{eq:asymsbm}
\overline{\boldsymbol{\pi}}_{\textnormal{SBM}}= \alpha \frac{1}{n} \ones + (1 - \alpha) \left( \mathbf{v} + \frac{\alpha \beta}{1 - \alpha \beta} (\mat v^T \mat u) \mathbf{u} \right) ,
\end{equation}
where $\beta \coloneqq \frac{p-q}{p+q},$ and $\mathbf{u} \in \mathbb{R}^n$ is a unit vector such that $u_i = \frac{1}{\sqrt{n}},$ for $i \in C_1$ and $u_i = -\frac{1}{\sqrt{n}}$ for $i \in C_2.$
\end{corollary}


\textbf{Proof:}
With equal-sized communities, i.e., $m = n/2$, we have $w_{\max} = w_{\min} = \frac{n}{2}(p+q).$ Therefore the conditions of Theorem \ref{theo:sbm} are satisfied if $p + q \gg \log^3(n)/n.$ Observe that the expected adjacency matrix can be written as $\overline{\mat{A}} = \frac{p+q}{2} \ones \ones^T + \frac{n}{2}(p-q) \mat u \mat u^T.$ Furthermore, $\mat W = \frac{n}{2}(p+q) \mat I.$ Therefore $\overline{\mat P} = \overline{\mat A} \mat W^{-1} = \frac{1}{n} \ones \ones^T + \frac{p-q}{p+q} \mat{u} \mat{u}^T.$ From (\ref{eq:prSBM}), the asymptotic PageRank $\overline{\boldsymbol{\pi}}_{\textnormal{sbm}}$ is therefore given as
\[
\overline{\boldsymbol{\pi}}_{\textnormal{sbm}} = \alpha \overline{\mat P} \overline{\boldsymbol{\pi}}_{\textnormal{sbm}} + (1- \alpha) \mat v.
\]
Consequently, $\overline{\boldsymbol{\pi}}_{\textnormal{sbm}} = \frac{\alpha}{n} \ones + \alpha \beta \mat u \mat u^T \overline{\boldsymbol{\pi}}_{\textnormal{sbm}} + (1 - \alpha) \mat v,$ or $\left[ \mat I - \alpha \beta \mat u \mat u^T\right] \overline{\boldsymbol{\pi}}_{\textnormal{sbm}} = \frac{\alpha}{n} \ones + (1 - \alpha) \mat v.$ By Woodbury Matrix Inversion Lemma in \cite{hornjohn12},  $\left[ \mat I  - \alpha \beta \mat u \mat u^T\right]^{-1} = \mat I + \frac{\alpha \beta}{1 - \alpha \beta} \mat u \mat u^T.$  Hence we obtain $\overline{\boldsymbol{\pi}}_{\textnormal{sbm}} = \frac{\alpha}{n} \ones + (1 - \alpha) \left (\mat v +  \frac{\alpha \beta}{1 - \alpha \beta} (\mat u^T \mat v) \mat u  \right),$ using the fact that $\mat u$ and $\ones$ are orthogonal vectors.
\qed
The above corollary asserts that on an SBM matrix the PageRank is well approximated in the asymptotic regime of large graph size by the convex combination of the uniform probability vector $\frac{1}{n} \ones$, which is the asymptotic stationary distribution of a simple random walk on the SBM graph, and a linear combination of the preference vector $\mat v$ and the projection of the preference vector onto the community partitioning vector $\mat u.$ Thus in this simple scenario of SBM graphs with equally sized communities, we observe that PageRank incorporates information about the community structure, in the form of a term correlated with the partition vector $\mat u,$ as opposed to the usual random walk, which misses this information. It can also be inferred from (\ref{eq:asymsbm}) that if the correlation between the preference vector $\mat v$ and $\mat u$ is large, e.g., when the seed set of PageRank is chosen to be in one of the communities, the resulting PageRank will display a clear delineation of the communities. This provides a mathematical rationale for why PageRank works for semi-supervised graph partitioning \cite{AGMS12}, at least in the asymptotic regime.

To prove Theorem \ref{theo:sbm} we need the following Lemmas, whose proofs are given in Appendix \ref{ap:secsbmprank}.

\begin{lemma}
\label{lemma:deg_conc}
For an SBM graph $\mathcal{G}(m,n-m,p,q),$ when $w_{\min} = \mathrm\omega(\log^3(n))$ it can be shown that for some $C,$
\[
\max_{ 1 \leq i \leq n} \left|\frac{D_i}{\mathbb{E} (D_i)} - 1 \right|  \leq C\sqrt{\frac{\log(n)}{ w_{\min}}} \text{ w.h.p}.
\]
\end{lemma}
The proof of this lemma follows from applying Bernstein's concentration lemma to the degrees of SBM graph. The proof is given in Appendix \ref{ap:lemsbmdeg}.

For ease of notation, let $\overline {\mathbf{Q}} = \mat W^{-1/2} \mathbb{E} (\mathbf{A}) \mat W^{-1/2},$ where $\mathbf{W} = \mathbb{E}(\mathbf{D}).$ As before $\mathbf{Q} = \mathbf{D}^{1/2} \mathbf{A} \mathbf{D}^{1/2}.$ We need the following concentration result on $\mat Q.$

\begin{lemma}\label{errorboundadj}
For an SBM graph for which $w_{\min} = \mathrm\omega(\log^3(n)),$ and $\frac{w_{\max}}{w_{\min}} \le C$ for some $C,$ it can be shown that
\[
\| \mathbf{Q} - \overline{\mathbf{Q}}\|_2 = C \frac{\sqrt{\log(n)w_{\max}}}{w_{\min}} = o(1)
\]
w.h.p.
\end{lemma}
We prove this lemma in Appendix \ref{ap:errorbadj}.\\
\textbf{Proof of Theorem \ref{theo:sbm}:}
We write the error between $\boldsymbol{\pi}$ and $\overline{\boldsymbol{\pi}}$ as follows
\begin{align}
\boldsymbol{\delta} &= \boldsymbol{\pi} - \overline{\boldsymbol{\pi}} \nonumber \\
&= (1 - \alpha) \left[\mat D^{1/2} (\mat I - \alpha \mat Q )^{-1} \mat D^{-1/2} - \mat W^{1/2} (\mat I - \alpha \overline{\mat Q})^{-1} \mat W^{-1/2} \right] \mat v \nonumber \\
&= (1-\alpha)\biggl[ \mat W^{1/2} \left((\mat I - \alpha {\mat Q} )^{-1} - (\mat I - \alpha \overline{\mat Q} )^{-1}\right) \mat W^{-1/2} \biggr ]\mat v +  \nonumber \\ &  (1-\alpha) \biggl[\mat D^{1/2} (\mat I - \alpha \mat Q )^{-1} \mat D^{-1/2} - \mat W^{1/2}(\mat I - \alpha \mat Q )^{-1} \mat W^{-1/2} \biggr ] \mat v , \label{eq:theosteps}
\end{align}
where in the last equality we added and subtracted $\mat W^{1/2}(\mat I - \alpha \mat Q )^{-1} \mat W^{-1/2}$ and reordered terms.
Now we analyse the two terms in square brackets in the last equality in (\ref{eq:theosteps}), which we denote $T_1$ and $T_2,$ respectively.
Notice that we have $\| \boldsymbol{\delta} \|_1 \le \|T_1\|_1 + \|T_2\|_1 .$ Next we show that as $n \to \infty,$ $\|T_1\|_1$ and $\|T_2\|_1$ are $o(1)$ separately and consequently we obtain the result of the theorem.

Let us first consider $T_1.$ We have
\begin{align*}
T_1 &=(1-\alpha)\biggl[ \mat W^{1/2} \left((\mat I - \alpha {\mat Q} )^{-1} - (\mat I - \alpha \overline{\mat Q} )^{-1}\right) \mat W^{-1/2} \biggr ]\mat v\\
&= (1- \alpha) \mathbf{W}^{1/2}(\mathbf{I} - \alpha {\mathbf{Q}})^{-1}\left( \overline{\mathbf{Q}} - \mathbf{Q} \right)(\mathbf{I} - \alpha \overline{\mathbf{Q}})^{-1}  \mathbf{W}^{-1/2} \mathbf{v},
\end{align*}
which we obtained by factoring out $(\mathbf{I} - \alpha {\mathbf{Q}})^{-1}$ and $(\mathbf{I} - \alpha \overline{\mathbf{Q}})^{-1}$ on  the left and right sides of the square brackets.
Next we focus on the 2-norm of $T_1.$

\begin{align*}
\|T_1\|_2 &\mylt{(a)} (1- \alpha) \sqrt{ w_{\max}} \|(\mathbf{I} - \alpha {\mathbf{Q}})^{-1}\|_2 \|\overline{\mathbf{Q}} - \mathbf{Q}\|_2 \|(\mathbf{I} - \alpha \overline{\mathbf{Q}})^{-1}\|_2 \frac{1}{\sqrt{w_{\min}}}\|\mathbf{v}\|_2 \\
&\mylt{(b)} \frac{1}{1-\alpha}\sqrt{\frac{w_{\max}}{w_{\min}}} \|\mat Q - \overline{\mat Q}\|_2\|\mat v\|_2 \\
&\mylt{(c)} C \frac{\sqrt{\log(n)w_{\max}}}{w_{\min}\sqrt{n}} \\
&= C \sqrt{\frac{\log(n)}{nw_{\max}}} \frac{w_{\max}}{w_{\min}}.
\end{align*}
This proves $\|T_1\|_1 \le \sqrt{n}\|T_1\|_2 \footnote{By Cauchy Schwartz inequality on norms.}\le C \sqrt{\frac{\log(n)}{w_{\max}}} \frac{w_{\max}}{w_{\min}} = o(1),$ from the assumptions of the theorem.
Here in (a) we used the submultiplicative property of matrix norms and the fact that 2-norm of diagonal matrices is the maximum diagonal element in magnitude. The inequality (b) follows because $\|(\mathbf{I} - \alpha {\mathbf{Q}})^{-1}\|_2 \le \frac{1}{1- \alpha}$ and $\|(\mathbf{I} - \alpha \overline{\mathbf{Q}})^{-1}\|_2 \le \frac{1}{1 - \alpha}$ and step (c) follows from Lemma \ref{errorboundadj} and the assumption that $\|\mat v\|_2 = O(1/\sqrt{n})$.

Next we analyse the second term $T_2.$ For ease of notation we denote $\widetilde{\mat R} = \mat W^{1/2} \left (  \mat{I} - \alpha \mat Q \right )^{-1} \mat W^{-1/2}.$ Then by simple algebraic manipulations
\begin{align*}
T_2 &= (1 - \alpha) \left [ \mat D^{1/2}  \left (  \mat{I} - \alpha \mat Q \right )^{-1} \mat D^{-1/2} - \mat W^{1/2} \left (  \mat{I} - \alpha \mat Q \right )^{-1} \mat W^{-1/2}\right]\mat v \\
&= (1-\alpha)\left(\mat D^{1/2} \mat W^{-1/2} \widetilde{\mat R} \mat W^{1/2} \mat D^{-1/2} - \widetilde{\mat R} \right) \mat v\\
&= (1-\alpha)\left(\mat D^{1/2} \mat W^{-1/2} \widetilde{\mat R} \left( \mat W^{1/2} \mat D^{-1/2} - \mat I \right) + \left( \mat D^{1/2}\mat W^{-1/2} - \mat I \right) \widetilde{\mat R} \right) \mat v,
\end{align*}
where the last step is obtained by adding and subtracting $\mat D^{1/2} \mat W^{-1/2} \widetilde{\mat R}.$

\noindent Now we have $\| \mat D^{1/2} \mat W^{-1/2} - \mat I \|_2 = \max_{i} \left |\sqrt{\frac{d_i}{w_i}} - 1 \right| \le \max_{i} \left |\frac{d_i}{w_i} - 1 \right| \le C \sqrt{\frac{\log(n)}{w_{\min}}} $ w.h.p. by Lemma \ref{lemma:deg_conc} and similarly
$\|\mat D^{1/2} \mat W^{-1/2} \|_2 \le \|\mat D^{1/2} \mat W^{-1/2} - \mat I\|_2 + \|\mat I\|_2 \le C \sqrt{\frac{\log(n)}{w_{\min}}} + 1.$ In addition $\|\mat W^{1/2} \mat{D}^{-1/2} - \mat I\|_2 = \max_i \left| \sqrt{\frac{w_i}{d_i}} - 1 \right| \le \max_i \left| \frac{w_i}{d_i} - 1 \right|.$ It can be shown that since $\max_i \left| \frac{d_i}{w_i} - 1 \right| \le C\sqrt{\frac{\log(n)}{w_{\min}}}$ w.h.p. (by Lemma \ref{lemma:deg_conc}), then $\max_i \left| \frac{w_i}{d_i} - 1 \right| \le C\sqrt{\frac{\log(n)}{w_{\min}}}$ w.h.p.\footnote{{\footnotesize This follows since we can write $\frac{d_i}{w_i} = 1 + \eta_i$, with $\max_i |\eta_i| = O\left(\sqrt{\frac{\log(n)}{w_{\min}}}\right) = o(1)$ w.h.p., then $\frac{w_i}{d_i} = \frac{1}{1+\eta_i} = 1 - \eta_i + O(\eta_i^2),$ hence $\max_i |\frac{w_i}{d_i} - 1| = O(\max_i |\eta_i|) = O\left(\sqrt{\frac{\log(n)}{w_{\min}}}\right) = o(1)$ w.h.p.}} Therefore $\|\mat W^{1/2} \mat{D}^{-1/2} \|_2 \le \|\mat W^{1/2} \mat{D}^{-1/2}  - \mat I \|_2 + \|\mat I\|_2 \le C\sqrt{\frac{\log(n)}{w_{\min}}} + 1 $ w.h.p.
Using the above facts and denoting $\delta = C\sqrt{\frac{\log(n)}{w_{\min}}} $ we obtain

\begin{align}
\|T_2\|_2 &\le \left (\|\mat D^{\frac{1}{2}}\mat W^{-\frac{1}{2}}\|_2 \|\widetilde{\mat R}\|_2 \|\mat W^{\frac{1}{2}} \mat D^{-\frac{1}{2}} - \mat I\|_2 + \|\mat D^{\frac{1}{2}} \mat W^{-\frac{1}{2}} - \mat I\|_2 \|\widetilde{\mat R}\|_2  \right)\|\mat v\|_2  \nonumber \\
\label{eq:proofstep}
&\le C(\delta(\delta + 1)\frac{1}{1-\alpha} + \delta) \frac{1}{1-\alpha}\sqrt{\frac{w_{\max}}{n w_{\min}}}\\
&\le C \delta\sqrt{\frac{w_{\max}}{n w_{\min}}} \mbox{w.h.p.}
\end{align}
Hence we have $\|T_2\|_1 \le \sqrt{n}\|T_2\|_2 \le C \delta\sqrt{\frac{w_{\max}}{w_{\min}}}$ w.h.p., which from our assumptions is $o(1).$
Here in (\ref{eq:proofstep}) we used the fact that
 \[
\|\widetilde{\mat R}\|_2 = \|\mat W^{1/2} \left (  \mat{I} - \alpha \mat Q \right )^{-1} \mat W^{-1/2}\|_2 \le \sqrt{\frac{w_{\max}}{w_{\min}}}\|\mat I  - \alpha \mat{Q}\|_2 \le \frac{1}{1 - \alpha} \sqrt{\frac{w_{\max}}{w_{\min}}} \le C,
\] and that $\|\mat v\|_2 \le C/\sqrt{n} ,$ for some $C.$
\qed
\vspace{5pt}
\noindent \textit{Remark:} This method of proof can be extended to similar models like the Stochastic Block Model with multiple communities and their generalizations, e.g., Random Dot Product Graphs \cite{athreya2013}.

\section{Experimental Results}

\begin{figure}
\centering
\includegraphics[scale=0.65]{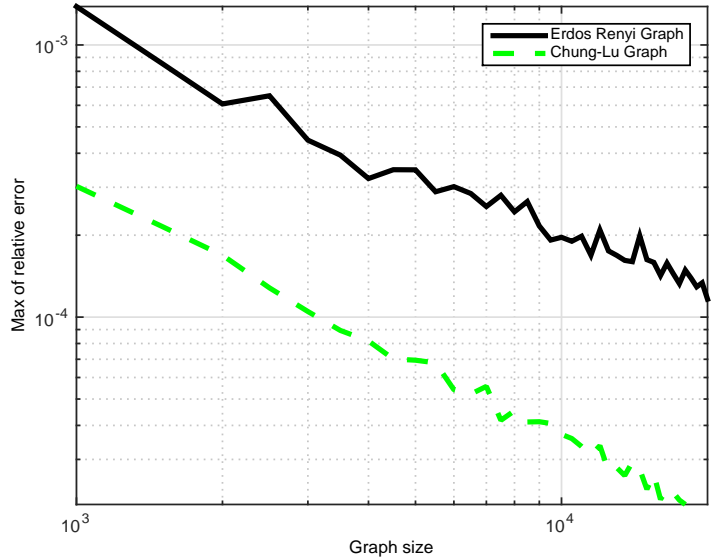}
\caption{Log-log plot of maximum normalized error for ER and Chung-Lu graphs}
\label{fig1}
\end{figure}

\begin{figure}
\centering
\includegraphics[scale=0.65]{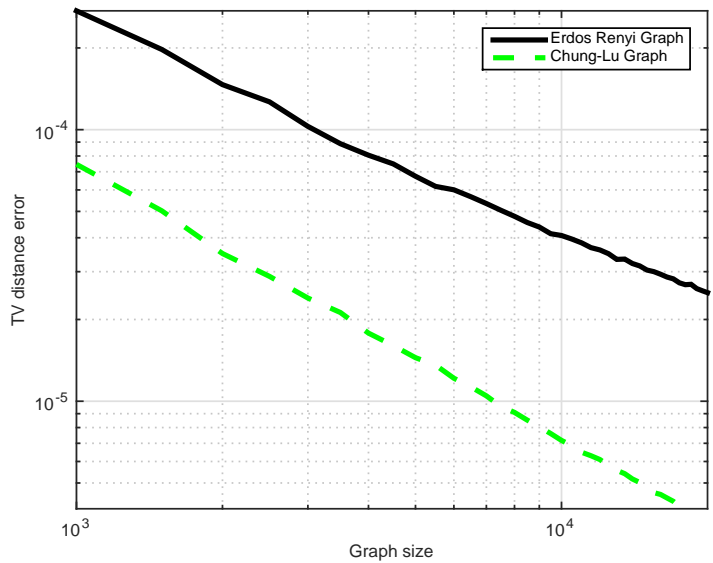}
\caption{Log-log plot of TV distance error for ER and Chung-Lu graphs}
\label{fig:tverror}
\end{figure}

In this section, we provide experimental evidence to further illustrate the analytic results obtained in the previous sections. In particular, we simulated
ER graphs with $p_n = C\frac{\log^7(n)}{n}$ and
Chung-Lu graphs with the degree vector $w$ sampled from a geometric distribution so that the average degree $\overline{w} = c n^{1/3},$ clipped such that $w_{\max}= 7w_{\min}$,
for various values of graph size, and plotted the maximum of normalized error $\widetilde{\delta}$ and TV distance error $\norm{\delta}_1$, respectively, in Figures \ref{fig1} and \ref{fig:tverror}. As expected, both these errors decay as functions of $n,$ which illustrates that the PageRank vector does converge to the asymptotic value.


\begin{figure}
\centering
\includegraphics[scale=0.5]{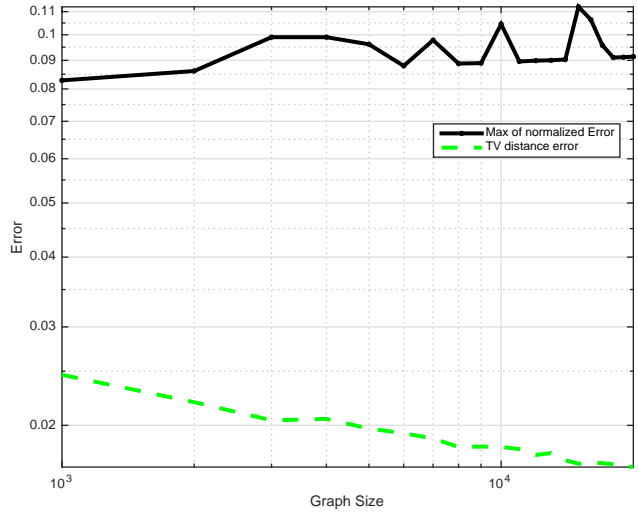}
\caption{Log-log plot of TV distance and maximum
error for power-law graphs}
\label{figpowererror}
\end{figure}

In the spirit of further exploration, we have also conducted simulations on power-law graphs with exponent $\beta = 4$ using the Chung-Lu graph model with $w_i = ci^{-1/(\beta - 1)},$ for $i_0 \le i \le n + i_0 $ with
\[
c = \frac{\beta - 2}{\beta - 1}dn^{1/(\beta - 1)},
\]
\[
i_0 = n \left[ \frac{d(\beta-1}{m(\beta-2)}\right]
\]
Please refer to \cite{CLV03} for details. We set max degree $m = n^{1/3}$ and average degree $d = n^{1/6}.$
In Figure \ref{figpowererror} we observe that for this graph the max-norm of the relative error does not converge to zero. On the other hand the TV-norm seems to converge to zero with graph size, albeit very slowly. Note that these graphs satisfy Property \ref{prop:fast_mixing} \cite{CLV03}, but they do not satisfy Property \ref{prop:bounded_degrees}. Therefore at this point, it is not possible to conclude whether the assumption of bounded variation of degrees is necessary for the convergence to hold. It might be interesting to investigate in detail the asymptotic behavior of  PageRank in undirected power-law graphs.

\begin{figure}
\centering
\includegraphics[scale=0.6]{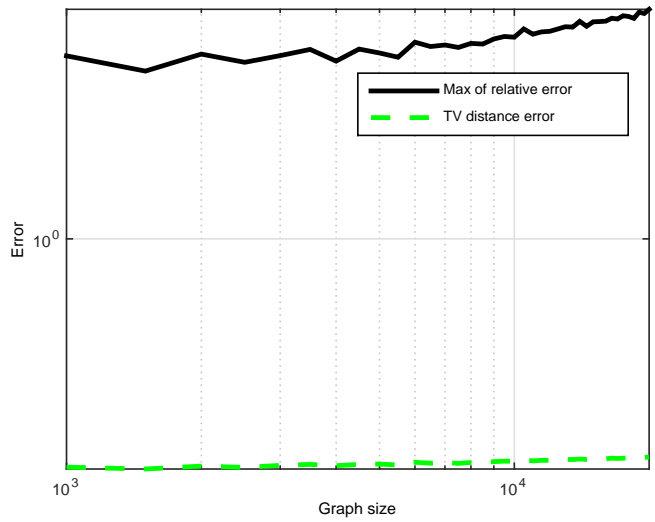}
\caption{Log-log plot of TV distance and maximum relative error for ER-graph when $v = e_1$}
\label{figv1}
\end{figure}

\begin{figure}
\centering
\includegraphics[scale=0.5]{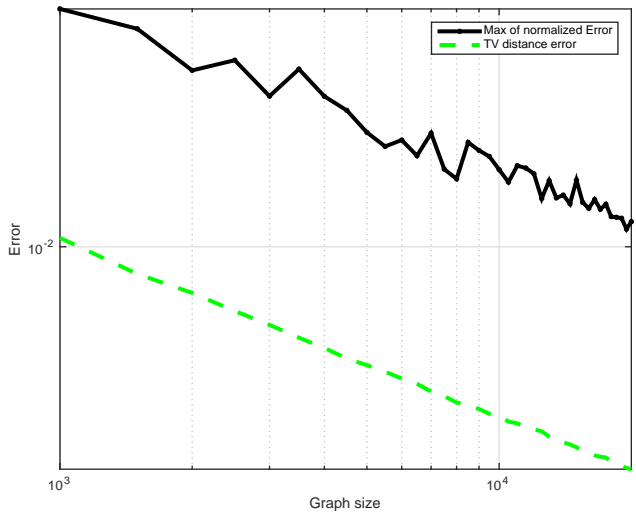}
\caption{Log-log plot of maximum normalized error and TV-distance error for an SBM graph}
\label{figsbm}
\end{figure}

Furthermore, we also see that in the case $\mat v= \mat e_i,$ the standard unit vector, for some $i$ we do not have the conjectured convergence, as can be seen on Figure \ref{figv1} in the case of ER graphs. It can also be seen from our analysis that if $v_k = 1 $ for some $k,$ the quantity $\norm{\widetilde{Q}D^{-1/2}v}_{\infty},$ becomes:
\vspace{-0.2cm}
\begin{equation*}
\max_i \left|\sum_j \left ( \frac{A_{ij}}{\sqrt{d_i d_j}} - \frac{\sqrt{d_i d_j}}{\sum_l d_l} \right ) v_j /\sqrt{d_j}\right|
= \max_{i} \frac{1}{\sqrt{d_i} d_k} \left|A_{ik} - \frac{d_i d_k}{\sum_l d_l}\right|,
\vspace{-0.1cm}
\end{equation*}
which is $O\left(\frac{1}{\sqrt{w_{\min}}w_{k}}\right)$ and does not fall sufficiently fast.
We simulated an SBM matrix with two communities of equal size, with $p = 0.1$ and $q = 0.01.$ In Figure \ref{figsbm} we plot the maximum normalized error and the TV-distance error against graph size on a log-log plot. As expected both errors go to zero for large graph sizes.

\section{Conclusions}
In this work, we have shown that when the size of a graph tends to infinity, the PageRank vector lends itself to be approximated by a mixture of the preference vector and the degree distribution, for a class of undirected random graphs including the Chung-Lu graph. We expect that these findings will shed more light on the behaviour of PageRank on undirected graphs, and possibly help to optimize the PageRank operation, or suggest further modifications to better capture local graph properties. We also obtain an asymptotic expression for the PageRank on SBM graphs. It is seen that this asymptotic expression contains information about community partitioning in the simple case of SBM with equal-sized communities. It would be interesting to study the implications of our results for community detection algorithms.
\vspace{-0.1cm}
\section*{Acknowledgements} We would like to thank Nelly Litvak for stimulating discussions on the topic of the paper. The work of K. Avrachenkov  and A. Kadavankandy was partly funded by the French Government (National Research Agency, ANR) through the ``Investments for the Future'' Program reference \#ANR-11-LABX-0031-01 and the work of L. Ostroumova Prokhorenkova and A. Raigorodskii was supported by Russian Science Foundation (\# 16-11-10014).

\appendix
\section{Proof of Lemma \ref{Bound on Sv}}\label{ap:lemmaSv}
From Lemma \ref{degree_concentration}, we have for Chung-Lu graphs that: $d_i = w_i (1+ \epsilon_i)$, where $\eta \equiv \max_i \epsilon_i = o(1)$ with high probability. In the proof we assume explicitly that $v_i = 1/n,$ but the results hold in the slightly more general case where $v_i = O(1/n)$ uniformly $\forall i$, i.e., $\exists K $ such that $\max_i n v_i \leq K.$ It can be verified easily that all the bounds that follow hold in this more general setting.
The event $\{ \eta = o(1)\},$ holds w.h.p. asymptotically from Lemma~\ref{degree_concentration}. In this case, we have
\begin{equation*}
\sum_j \left ( \frac{A_{ij}}{\sqrt{d_i d_j}}  - \frac{\sqrt{d_i d_j}}{\sum_i d_i} \right) \frac{v_j}{\sqrt{d_j}} =
\sum_j \left ( \frac{A_{ij}}{\sqrt{d_i d_j}}  - \frac{\sqrt{d_i d_j}}{\sum_k d_k} \right) \frac{v_j}{\sqrt{w_j}}(1 + \varepsilon_{j})
\end{equation*}
where $\varepsilon_{j}$ is the error of convergence, and we have $\max_{j} \varepsilon_{j} = O(\eta)$. Therefore,
\begin{align}
\norm{\widetilde{\mat Q}\mat v'}_{\infty} &\leq \norm{\widetilde{\mat Q} \mat q}_{\infty} + \max_i \varepsilon_i  \norm{\widetilde{\mat Q} \mat q}_{\infty} \nonumber \\
&\leq  \norm{\widetilde{\mat Q} \mat q}_{\infty}(1 + o(1)) \quad \mbox{w.h.p.}, \label{vequival}
\end{align}
where $\mat q$ is a vector such that $q_i = \frac{nv_i}{\sqrt{w_i}}.$
Furthermore, we have w.h.p.
\begin{align*}
\frac{A_{ij}}{\sqrt{d_i d_j}}  - \frac{\sqrt{d_i d_j}}{\sum_k d_k} & = \frac{A_{ij}}{\sqrt{w_i(1 + \epsilon_i) w_j(1 + \epsilon_j)}}  - \frac{\sqrt{w_i(1 + \epsilon_i) w_j(1 + \epsilon_j)}}{\sum_k w_k(1 + \epsilon_k)} \\
&= \frac{A_{ij}}{\sqrt{w_i w_j}} \left( 1 + O(\epsilon_{i}) + O(\epsilon_{j}) \right) - \frac{\sqrt{w_i w_j}}{\sum_k w_k} \left (\frac{ 1 + O(\epsilon_{i}) + O(\epsilon_{j})}{1 + O(\eta)}\right) \\
& = \left ( \frac{A_{ij}}{\sqrt{w_i w_j}} - \frac{\sqrt{w_i w_j}}{\sum_k w_k} \right)(1 + \delta_{ij}),
\end{align*}
where $\delta_{ij}$ is the error in the $ij^{\text{th}}$ term of the matrix and $\delta_{ij} = O(\eta)$ uniformly, so that $\max_{ij} \delta_{ij} = o(1)$ w.h.p.
Consequently, defining $\widetilde{\overline{Q}}_{ij} = \frac{A_{ij}}{\sqrt{w_i w_j}} - \frac{\sqrt{w_i w_j}}{\sum_k w_k} $ we have:
\begin{align}
\norm{\widetilde{\mat Q} \mat q}_{\infty} &\leq \norm{\widetilde{\overline{\mat Q}} \mat q}_{\infty} + \max_i |\sum_j \widetilde{\overline{Q}}_{ij} \delta_{ij} q_j| \nonumber  \\
&\leq \norm{\widetilde{\overline{\mat Q}} \mat q}_{\infty} +  O(\eta) \max_i \frac{1}{\sqrt{w_{\min}}} \sum_j |\widetilde{\overline{\mat Q}}_{ij}|  \nonumber  \\
& \leq \norm{\widetilde{\overline{\mat Q}} \mat q}_{\infty}  + o(1) \frac{1}{\sqrt{w_{\min}}}  \left(C\sqrt{\frac{w_{\max}}{w_{\min}}} + \frac{w_{\max}}{w_{\min}}\right)\label{infbound} \\
& \leq \norm{\widetilde{\overline{\mat Q}} \mat q}_{\infty}  + o(1/\sqrt{w_{\min}}) \label{unibound}
\end{align}
where in (\ref{infbound}) we used the fact the $O(\eta)$ is a uniform bound on the error and it is $o(1)$ w.h.p. and $\max_j q_j \leq \frac{1}{\sqrt{w_{\min}}}.$ In (\ref{infbound}) we also used the fact that
\begin{align*}
\max_{i}\sum_j|\widetilde{\overline{Q}}_{ij}|
&\leq \max_i \sum_j \frac{A_{ij}}{\sqrt{w_i w_j}} + \sum_j \frac{\sqrt{w_i w_j}}{\sum_k w_k} \\
&\leq \max_i \frac{1}{\sqrt{w_{\min}}} \frac{d_i}{\sqrt{w_i}} + \max_i \frac{\sqrt{w_i w_{\max}}}{w_{\min}}\\
&\mylt{(a)} C\sqrt{\frac{w_i}{w_{\min}}} + \frac{w_{\max}}{w_{\min}}\\
&\leq C\sqrt{\frac{w_{\max}}{w_{\min}}} + \frac{w_{\max}}{w_{\min}},
\end{align*}
where $C$ is some constant. In (a) above we used the fact that w.h.p. $d_i = w_i (1 + o(1)),$ by Lemma \ref{degree_concentration}, hence $\exists C$ such that $\forall n$ large enough $d_i \le C w_i.$

Now we proceed to bound $\norm{\widetilde{\overline{\mat Q}} \mat q}_{\infty}.$
Substituting for $q_i = \frac{1}{\sqrt{w_i}},$ we get
\begin{align}
\sum_j \frac{1}{\sqrt{w_j}} \left( \frac{A_{ij}}{\sqrt{w_i w_j}} - \frac{\sqrt{w_i w_j}}{\sum_k w_k} \right)
= \sum_j \frac{1}{w_j \sqrt{w_i}} \left ( A_{ij} - \frac{w_i w_j}{\sum_i w_i} \right) \nonumber \\
 \equiv \frac{1}{\sqrt{w_i}}X_i \label{final_reln}.
\end{align}
We seek to bound $\max_i |X_i|$:
\[
X_i = \sum_j \frac{1}{w_j} \left ( A_{ij} - \frac{w_i w_j}{\sum_i w_i} \right).
\]

Furthermore, $\mathbb{E}(X_i^2) = \sum_j \frac{1}{w_j^2} \mathbb{E}(A_{ij} - p_{ij})^2,$ with $p_{ij} = \frac{w_i w_j}{\sum w_i}.$
So, $\mathbb{E}(X_i^2) = \sum_j \frac{1}{w_j^2} p_{ij}(1-p_{ij}) \leq \frac{w_i}{\sum_i w_i} \sum_j \frac{1}{w_j} \leq n \frac{p_i}{w_{\min}},$ where $p_i = \frac{w_i}{\sum_i w_i},$ and $\frac{A_{ij}}{w_j} \leq 1/w_{\min}.$ Therefore using Bernstein Concentration Lemma for $\epsilon < n\max_i p_i$:
\begin{align}
\mathbb{P} \left ( \max_i |\sum_j (A_{ij} - p_{ij})/w_j| \geq \epsilon \right ) & \leq n \max_i \exp(-\frac{\epsilon^2}{2(p_in /w_{\min}) + \epsilon/w_{\min}})\nonumber \\
     & \leq n \max_i \exp(-\frac{w_{\min} \epsilon^2}{ 2(np_i + \epsilon)}) \nonumber \\
     & \leq  n  \exp(- \epsilon^2 w_{\min}/(4n\max_i p_i)) \nonumber \\
     & \leq n \exp(\frac{-\epsilon^2 \text{vol} w_{\min}}{4w_{\max}n}), \label{eq:finstep}
\end{align}
where $\frac{\text{vol}}{n} = \frac{\sum_i w_i}{n} \geq w_{\min}.$
It can be verified that when $\epsilon = \frac{1}{(\overline{w})^\alpha}$, for some $\alpha > 0,$ the RHS of (\ref{eq:finstep}) can be upper bounded by $n^{-(\gamma K - 1)},$ if $\overline{w} \geq (\gamma \log(n))^{\frac{1}{1 - 2\alpha}},$ for some large enough $\gamma$, which can be easily satisfied if $w_{\min} \gg O(\log^c(n)),$ for some $c > 1,$ where $K$ is a constant such that $w_{\max} \le K w_{\min}.$
Thus, finally, from (\ref{final_reln}) and (\ref{unibound}) we have $\norm{\widetilde{\mat Q}\mat q}_{\infty} = o(1/\sqrt{w_{\min}}),$ w.h.p., and therefore from (\ref{vequival}), we get $\norm{\widetilde{\mat Q}\mat v'}_{\infty} = o(1/\sqrt{w_{\min}}).$

\qed

\section{Proof of Lemmas in Section \ref{sec:sbmprank}}\label{ap:secsbmprank}
\subsection{Proof of Lemma \ref{lemma:deg_conc}}\label{ap:lemsbmdeg}
 The proof is an application of Bernstein's Concentration Lemma. Note that for $1 \leq i \leq m,$ $D_i = \sum_j A_{ij}$ . Here the mean degree $ \mathbb{E} (D_i) = m p + (n-m)q = t_1,$ and the variance $B^2 _n = mp(1-p) + (n-m)q(1-q) \le t_1$ for $i \le m.$ Similarly for $i >m,$ $\mathbb{E}(D_i) = m q + (n-m)p = t_2$ is  and variance $\mbox{Var}[{D_i}] \le t_2.$ Then, the minimum average degree $w_{\min} = \min(t_1,t_2).$
By Bernstein's Lemma, for $\epsilon = C\sqrt{\frac{\log(n)}{ w_{\min}}},$
\begin{align*}
\mathbb{P} \left(\max_{1 \le i \le m} |D_i - t_1| \ge \epsilon t_1 \right) &\le 2 m \exp\left(\frac{-\epsilon^2 t_1^2}{2(t_1 \epsilon/3+ t_1)}\right) \\
& = 2 m \exp\left(\frac{-\epsilon^2 t_1}{1 + \epsilon/3}\right)\\
& = O(n^{-c}),
\end{align*}
for some $c.$ Hence $\max_{1 \le i \le m} \left | \frac{D_i - t_1}{t_1} \right| \le C\sqrt{\frac{\log(n)}{w_{\min}}}$ w.h.p.
Similarly
\[
\max_{1+m \le i \le n/2} \left|\frac{D_i - t_2}{t_2} \right| \le C\sqrt{\frac{\log(n)}{w_{\min}}}, \mbox{ w.h.p}.
\]
Combining the two bounds above we get,
\begin{equation}
\label{eq:degree_concentration}
 \max_{ 1 \leq i \leq n} \left|\frac{D_i}{\mathbb{E} (D_i)} - 1 \right|  \leq C\sqrt{\frac{\log(n)}{ w_{\min}}}, \text{ w.h.p}.
\end{equation}
\qed

\subsection{Proof of Lemma \ref{errorboundadj}}\label{ap:errorbadj}
To prove Lemma \ref{errorboundadj} we need the following lemma on the spectral norm of the difference between the adjacency matrix and its mean.

\begin{lemma}
\label{lem:specbound}
For an SBM matrix $G(m,n-m,p,q)$  with adjacency matrix $\mat A$ and $\overline{\mat A} = \mathbb{E} (\mat A),$ there exists a constant $K$ s.t.
\[
\|\mat A - \overline{\mat A}\|_2 \le K\sqrt{\log(n) w_{max}}, \mbox{  w.h.p.},
\]
where $w_{max} = \max(m,n-m)p + \min(m,n-m)q$ is the maximum average degree,
if $w_{\max} = \mathrm\omega(\log^3(n)).$
\end{lemma}

To prove this Lemma we need the Matrix Bernstein Concentration result, which we state below for the sake of completeness:
\begin{lemma}
\label{lem:matbernstein}
\cite[Theorem ~1.4]{tropp2012user}.
Let $\mathbf{S}_1,\mathbf{S}_2,\ldots \mathbf{S}_t$ be independent random matrices with common dimension $d_1 \times d_2.$ Assume that each matrix has bounded deviation from its mean, i.e.,
$$\| \mathbf{S}_k - \mathbb{E} (\mathbf{S}_k) \| \le R, \text{ for each } k = 1,\ldots n.$$
Let $\mathbf{Z} = \sum_{k=1}^t \mathbf{S}_k $ and introduce a variance parameter
\begin{equation*}
\sigma^2_{\mathbf{Z}} = \max \left \{ \|\mathbb{E}\left( (\mathbf{Z} - \mathbb{E} (\mathbf{Z})) (\mathbf{Z} - \mathbb{E} (\mathbf{Z}))^H\right) \|, \|\mathbb{E}\left( (\mathbf{Z} - \mathbb{E} (\mathbf{Z}))^H (\mathbf{Z} - \mathbb{E} (\mathbf{Z}))\right) \| \right \}.
\end{equation*}
Then
\begin{equation}
\label{eq:matbernstein}
\mathbb{P} \{ \| \mathbf{Z} - \mathbb{E} (\mathbf{Z}) \| > t \} \le (d_1 + d_2) .\exp\left ( \frac{-t^2/2}{\sigma^2_{\mathbf{Z}} + R t/3} \right),
\end{equation}
$\text{ for all } t \ge 0.$
\end{lemma}
\textbf{Proof of Lemma \ref{lem:specbound}:}
With $\mathbf{Z} = \mat A,$ in Lemma \ref{lem:matbernstein}, we can decompose $\mathbf{Z}$ as sums of Hermitian matrices $\mathbf{S}_{i^{'}j^{'}},$ $\mathbf{Z} = \sum_{1 \le i^{'}<j^{'} \le n} \mat{S}_{i^{'}j^{'}}$ such that:
\begin{equation}
(\mathbf{S}_{i^{'}j^{'}})_{ij}
=
\begin{cases}
 A_{i^{'}j^{'}} \text{ if } i = i^{'}, j = j^{'}, \\
 A_{i^{'}j^{'}} \text{ if } i = j^{'}, j = i^{'}, \\
0 \text{ otherwise.}
\end{cases}
\end{equation}
Notice that if $\mat x \neq 0,$ $\|(\mathbf{S}_{i^{'}j^{'}} - \mathbb{E} (\mathbf{S}_{i^{'}j^{'}}) ) \mathbf{x} \|_2 = |2x_{i^{'}} x_{j^{'}} (A_{i^{'}j^{'}} - \mathbb{E}(A_{i^{'}j^{'}}))| < | x_{i^{'}}^2 + x_{j^{'}}^2 |.$ Consequently $\| \mathbf{S}_{i^{'}j^{'}} - \mathbb{E} (\mathbf{S}_{i^{'}j^{'}}) \|_2 < 1,$ giving $R = 1$ in the statement of Lemma \ref{eq:matbernstein}. Let $\mat{Y} = \mathbb{E}\left( (\mathbf{Z} - \mathbb{E} \mathbf{Z})^H (\mathbf{Z} - \mathbb{E} \mathbf{Z})\right),$ then
\begin{equation}
Y_{ij} =
\begin{cases}
v_1  & \text{ if } i = j  , i \le m,\\
v_2  & \text{ if }  i = j, i > m, \\
0  & \text{ otherwise},
\end{cases}
\end{equation}
where $v_1 =  m p(1-p)+ q(1-q)(n - m), v_2 = (n-m)p(1-p)+m q(1-q) .$
Therefore $\sigma^2_{\mathbf{Z}} = \max(v_1,v_2) = \max(n-m,m) p + \min(n-m,m) q =  \sigma^2.$ By our assumptions on the probabilities, $\sigma^2 = \mathrm\omega(\log^3(n)).$ Thus it follows that
\begin{align*}
\mathbb{P} (\| \mat A - \overline{\mat A} \| \geq t \sigma)  &\leq  2n \exp\left(\frac{-t^2 \sigma^2}{ 2\sigma^2 + t \sigma/3}\right) \\
&\leq 2n \exp(-t^2/3),
\end{align*}
if $\sigma > t.$ The RHS is $O(n^{-c})$ if $t >  \sqrt{r\log(n)},$ for some $r.$ \qed
Finally we are in a position to prove Lemma \ref{errorboundadj}\\
\textbf{Proof of Lemma \ref{errorboundadj}:}
We prove this result in two steps.
First we show that
\begin{equation}
\label{eq:step1}
\|\mat D^{-1/2} \mat A \mat D^{-1/2} - \mat W^{-1/2} \mat A \mat W^{-1/2} \|_2 = C \sqrt{\frac{\log(n)}{ w_{\min}}} = o(1).
\end{equation}
Observe that
\begin{align*}
&\|\mat D^{-1/2} \mat A \mat D^{-1/2} - \mat W^{-1/2} \mat A \mat W^{-1/2} \|_2 = \|\mat{Q} - \mat W^{-1/2} \mat D^{1/2} \mat Q \mat D^{1/2} \mat W^{-1/2} \| \\
& = \| \mat Q -  \mat W^{-1/2} \mat D^{1/2} \mat Q + \mat W^{-1/2} \mat D^{1/2} \mat Q - \mat W^{-1/2} \mat D^{1/2} \mat Q \mat D^{1/2} \mat W^{-1/2} \|_2 \\
& = \|(\mat I - \mat W^{-1/2} \mat D^{1/2}) \mat Q + \mat W^{-1/2}\mat D^{1/2} \mat Q(\mat I - \mat D^{1/2} \mat W^{-1/2}) \|_2\\
& \le \delta + (1 + \delta) \delta,
\end{align*}
where $\delta = \max_i \left| \frac{d_i}{w_i} - 1\right|.$
In the last line we used the fact that $\|\mat Q\|_2 = 1, \|\mat I - \mat W^{-1/2} \mat D^{1/2}\|_2 = \max_i \left| \sqrt{\frac{d_i}{w_i}} - 1\right| \le \max_i \left| \frac{d_i}{w_i}- 1\right|$ and

\[
\| \mat W^{-1/2}\mat D^{1/2}\|_2 \le \|\mat W^{-1/2}\mat D^{1/2} - \mat I\|_2 + \|\mat I\|_2 \le \delta + 1.
\]

By Lemma \ref{lemma:deg_conc}, $\delta \le  C\sqrt{\frac{\log(n)}{ w_{\min}}} = o(1)$ w.h.p.
Next we show that

\begin{equation}
\label{eq:step2}
\|\mat W^{-1/2} \mat A \mat W ^{-1/2} - \mat W^{-1/2} \overline{\mat A} \mat W^{-1/2}\|_2 \le \frac{C \sqrt{\log(n) w_{\max}}}{w_{\min}} = o(1).
\end{equation}

Now using Lemma \ref{lem:specbound} we have
\begin{align*}
\|\mat W^{-1/2}\mat A \mat W^{-1/2} - \mat W^{-1/2} \overline{\mat A} \mat W^{-1/2} \| &\le \frac{\|\mat A  - \overline{\mat A} \|_2}{w_{\min}} \\
&\le \frac{c\sqrt{\log(n) w_{\max}}}{w_{\min}} \\
&= o(1), \mbox{ w.h.p.},
\end{align*}
if $w_{\min} = \mathrm{\omega}(\sqrt{\log(n) w_{\max}}),$ which is satisfied when $w_{\max} \le C w_{\min}$ for some $C,$ and $w_{\max} = \mathrm\omega(\log^3(n)).$
The result of Lemma \ref{errorboundadj} then follows from (\ref{eq:step1}) and (\ref{eq:step2}) by applying the triangular inequality.\qed
\printnomenclature[10em]
\end{document}